\documentclass[11pt,letterpaper,twoside,headrule]{article}
\usepackage{document}
\usepackage{amssymb,amsmath,amscd,amsthm,pb-diagram}

\newcommand{\HH}{{\mathbb H}}

\newcommand{\R}{{\mathbb R}}
\newcommand{\Z}{{\mathbb Z}}

\DeclareMathOperator{\wh}{Wh}
\DeclareMathOperator{\colim}{colim}
\DeclareMathOperator{\met}{Met}
\DeclareMathOperator{\sectional}{sec}
\DeclareMathOperator{\isom}{Isom}
\DeclareMathOperator{\vc}{vc}

\theoremstyle{definition}

\newtheorem*{BC}{Borel Conjecture for $N$}

\begin{document}


\def\evenhead{{\protect\centerline{\textsl{\large{James F. Davis}}}\hfill}}

\def\oddhead{{\protect\centerline{\textsl{\large{The Work of Tom Farrell and Lowell Jones in Topology and Geometry}}}\hfill}}

\pagestyle{myheadings} \markboth{\evenhead}{\oddhead}

\thispagestyle{empty} \noindent{{\small\rm Pure and Applied
Mathematics Quarterly\\ Volume 8, Number 1\\ (\textit{Special Issue:
In honor of \\ F. Thomas Farrell and Lowell E. Jones, Part 1 of 2})\\
1---14, 2012} \vspace*{1.5cm} \normalsize

\begin{center}
\Large{\bf The Work of Tom Farrell and Lowell Jones in Topology and
Geometry}
\end{center}

\renewcommand{\thefootnote}{\fnsymbol{footnote}}

\begin{center}
{\large James F. Davis\footnote{I thank Shmuel Weinberger and Tom
Farrell for their helpful comments on a draft of this article.  I
also profited from reading the survey articles \cite{FJ(1990)} and
\cite{F(2008)}.}}
\end{center}


Tom Farrell and Lowell Jones caused a paradigm shift in high-dimensional topology, away from the view that high-dimensional topology was, at its core, an algebraic subject, to the current view that geometry, dynamics, and analysis, as well as algebra, are
key for classifying manifolds whose  fundamental group is infinite.
Their collaboration produced about fifty papers over a twenty-five year period.  In this tribute for the special issue of Pure and Applied Mathematics Quarterly in their honor, I will survey {\em some} of the impact of their joint work and mention briefly their individual contributions -- they have written about one hundred non-joint papers.

\section{Setting the stage}

In order to indicate the Farrell--Jones shift, it is necessary to describe the situation before the onset of their collaboration.  This is intimidating -- during the period of twenty-five years starting in the early fifties, manifold theory was perhaps the most active and dynamic area of mathematics.  Any narrative will have omissions and be non-linear.  Manifold theory deals with the classification of manifolds.  There is an existence question -- when is there a closed manifold within a particular homotopy type, and a uniqueness question, what is the classification of manifolds within a homotopy type?

The fifties were the foundational decade of manifold theory.
Ren\'e Thom developed transversality, the basic tool of differential topology.  John Milnor showed there were exotic differentiable structures on the 7-sphere.  After the work of Milnor it became clear that one must distinguish between three different types of manifolds: smooth, topological, and piecewise linear (PL) and as well as three sorts of classification: up to diffeomorphism, homeomorphism, and PL-homeomorphism.

An important theme in the study of infinite groups is the  phenomena of geometric rigidity of discrete subgroups of Lie groups: given two abstractly isomorphic discrete subgroups of a Lie group, under what circumstances are they necessarily conjugate?  There is a connection between rigidity and aspherical manifolds since, according to Kenkichi Iwasawa, a simply-connected real Lie group $G$ has a maximal compact subgroup $K$, and the homogenous space $G/K$ is diffeomorphic to euclidean space.  Hence if $\Gamma$ is a discrete, cocompact, torsion-free subgroup of $G$, the space $\Gamma\backslash G /K$ is a closed, aspherical, smooth manifold, where aspherical means that the universal cover is contractible.  It is not difficult to show that two aspherical manifolds with isomorphic fundamental groups are homotopy equivalent.

In 1954 George Mostow showed that if  $G$ is a simply-connected solvable Lie group, then any two discrete, cocompact, torsion-free, abstractly isomorphic subgroups of $G$ have diffeomorphic double coset manifolds, but need not be conjugate.

Armand Borel, seeing the impossibility of making a general group theoretic conjecture, made the far-reaching topological conjecture.

\begin{BC}
Let $M$ and $N$ be  closed, aspherical topological manifolds.  Then any homotopy equivalence $h: M \to N$ is homotopic to a homeomorphism.
\end{BC}

Borel had little evidence for his topological rigidity conjecture.  But, now, largely due to the work of Farrell and Jones, it is now known for most manifolds $N$ of geometric interest.

The decade of the sixties was the golden age of high-dimensional manifold theory.  It is perhaps no accident that Farrell and Jones were students during this time.  Steven Smale developed handlebody techniques, and thereby proved the high-dimensional Poincar\'e conjecture and the $h$-cobordism theorem in dimensions greater than four, indicating, following Borel, that homotopy equivalence and homeomorphism are not that far apart.  Michel Kervaire and John Milnor developed the foundations of surgery theory, as applied to the classification of exotic spheres.  They showed, in particular, that in each dimension greater than four, there are a finite number of smooth structures on the sphere.  William Browder and Sergei Novikov made surgery theory into a systematic theory for classifying simply-connected manifolds.  Dennis Sullivan, motivated by the example of manifolds homotopy equivalent to complex projective space, set up homotopy theoretic foundations for the theory of Browder and Novikov.  There were three results in the sixties which gave early indications of non-simply connected techniques.  First, in the non-simply-connected case, obstructions to a $h$-cobordism being a product land in the Whitehead group $\wh(\pi) := K_1(\Z \pi)/\langle \pm g \in \pi \rangle$.  Second, the Fields medal result of Novikov that rational Pontryagin classes are homeomorphism invariant required an analysis of the torus and led Novikov to his famous conjecture stipulating the precise extent to which the rational Pontryagin classes are invariant under a homotopy equivalence of manifolds.   Third, the thesis of Farrell \cite{F(1971)} studied when a manifold is a fiber bundle over a circle.  Although his thesis did not use surgery theory explicitly, it became a foundational tool for studying non-simply-connected manifolds.

At the end of the sixties, there were two explosive developments.
Rob Kirby and Larry Siebenmann showed that transversality and
handlebody theory, and thereby their myriad consequences, worked in
the topological category and C.T.C. Wall wrote his magnificent tome
{\em Surgery on compact manifolds}, the foundational work on the
surgery classification of non-simply connected manifolds.

In summary, the basic outline of the surgery theoretic classification of high-dimensional manifolds was completed by the end of the decade; classification was reduced to homotopy theory and obstructions lying in the algebraic $K$- and $L$-groups $K_*(\Z \pi)$ and $L_*(\Z \pi)$ where $\pi$ is the fundamental group of the manifold.  However,  little was known about these $K$- and $L$-groups outside of the simply-connected case.  $K$-theory measures the extent to which a homotopy equivalence is a simple homotopy equivalence and $L$-theory measures the extent to which a simple homotopy equivalence is homotopic to a homeomorphism.  More precisely, if the Whitehead group $\wh(\pi)$ vanishes then every homotopy equivalence is a simple homotopy equivalence and every $h$-cobordism is a product, and if the assembly map (defined by Frank Quinn and Andrew Ranicki) $H_*(B\pi; \bf L) \to L_*(\Z\pi)$ is an isomorphism then every simple homotopy equivalence between closed aspherical manifolds with fundamental group $\pi$ is homotopic to a homeomorphism.  Both of these conditions, the vanishing of the Whitehead group and the isomorphism of the $L$-theory assembly map, are conjectured to be true for fundamental groups of closed aspherical manifolds, and, more generally, for torsion-free groups.

The first non-simply-connected case to be understood was  $\Z^n$, and thereby the Borel conjecture for the torus was proven. The $K$-theory was analyzed by Hyman Bass, Alex Heller, and Richard Swan, who showed that $\wh(\Z^n) = 0$.  The $L$-theory of $\Z^n$ was computed independently by Julius Shaneson and  Wall using Farrell's thesis and by Novikov using algebraic techniques.    Wall outlined a program for computation of the algebraic $K$- and $L$-groups of integral group rings of finite groups, or at least for reducing the computations to computations in algebraic number theory.  Friedhelm Waldhausen and Sylvain Cappell computed the $K$- and $L$-theory of amalgamated products and HNN extensions in terms of the algebraic $K$- and $L$-theory of the constituent pieces and the Nil and UNil-groups, which reflect nilpotent phenomena.  Farrell, both alone and in collaboration with Wu-Chung Hsiang was one of the first to examine these still mysterious Nil groups \cite{FH(1970)}, \cite{FH(1973)}, \cite{F(1977)}, \cite{F(1979)}.

Of the many applications of surgery theory to the classification of manifolds and to transformation groups, I will only mention the program of Jones \cite{J(1971)} to characterize the fixed point sets of cyclic groups on disks and spheres.

\section{Geometry in service of topology}

Before I move to the work of Farrell and Jones, I need to mention controlled topology.  I will only state the controlled $h$-cobordism theorem of Steve Ferry and not the closely related results of Tom Chapman and Frank Quinn.  If $(Y,d)$ is a metric space and $\varepsilon > 0$, then a homotopy $H : X \times I \to Y$ is {\em  $\varepsilon$-controlled} if for all $x \in X$, the tracks $H(\{x\} \times I)$ of the homotopy all have diameter less than $\varepsilon$ (measured in $Y$).  If $N_-$ is a Riemannian manifold, an $h$-cobordism $(W; N_-, N_+)$ is {\em $\varepsilon$-controlled} if the composite of the deformation retract $W \times I \to W$ of $W$ to $N_-$ with the retract $W \to N_-$ is $\varepsilon$-controlled.  The controlled $h$-cobordism theorem then says that when the dimension of $N_-$ is greater than four there is an $\varepsilon = \varepsilon(N_-) > 0$ so that every $\varepsilon$-controlled $h$-cobordism is a product.

We need to discuss two more developments before I finally come to the paradigm shift of Farrell and Jones.  The first was that Farrell and Jones started collaborating (it was natural for them to do so, after all they were academic brothers -- both students of Hsiang).  They collaborated on a variety of problems in dynamics.  Farrell and Jones \cite{FJ(1978)} constructed an Anosov diffeomorphism on a connected sum of a torus with an exotic sphere, providing the first example of an Anosov diffeomorphism on a manifold which was not infranil.  The second was Farrell and Hsiang's proof \cite{FH(1983)flat} of the Borel conjecture for flat and almost flat manifolds of dimension greater than four.  This proof included the use of expanding endomorphisms and controlled topology to show that the Whitehead group vanishes, hinting at the possibility of applications of geometry to topology.

So what was the paradigm shift of Farrell and Jones?  Well, prior to their joint work, techniques for classifying high-dimensional manifolds included algebra, number theory, algebraic topology, and geometric topology.  Farrell and Jones introduced ideas from dynamics and differential geometry and applied them to questions such as the Borel Conjecture.\footnote{Analogous paradigm shifts occurred in the topology of low-dimensional manifolds.  In particular, William Thurston introduced geometric ideas to the classification of three-manifolds and Simon Donaldson introduced gauge theoretic ideas motivated by physics and applied them to the smooth classification of four-manifolds.}  It is notable that the tools they introduced to the study of topological rigidity are some of the key tools for geometric rigidity.

Their first amazing result was the verification of the Borel Conjecture for hyperbolic manifolds of dimensions greater than four \cite{FJ(1986)} and \cite{FJ(1989)Borel}.  I will outline this, not so much to give a coherent account, but rather to indicate the scope of ideas.  Let $M$ be a closed hyperbolic manifold.  The key geometric idea of Farrell and Jones was the notion of the asymptotic transfer:   given a path $\alpha : I \to M$ and a unit tangent vector $v \in S_{\alpha(0)}M$, the asymptotic transfer is a unit tangent vector field along $\alpha$ given by a map $v\alpha : I \to SM$.  This asymptotic transfer is not the vector field given by parallel translation but rather a clever variant defined using the identification of the universal cover $\widetilde M$ with hyperbolic space $\HH$.  Lift both the path $\alpha$ and the initial vector $v$ to hyperbolic space: $\widetilde \alpha : I \to \HH$ and $\widetilde v \in S_{\widetilde \alpha(0)}\HH$.  Let $\gamma: [0,\infty] \to \HH \cup S_\infty$ be the geodesic ray with initial tangent vector $\widetilde v$.  Then define $\widetilde{v\alpha} : I \to S\HH$ so that $\widetilde{v\alpha}(t) \in S_{\widetilde \alpha(t)}\HH$ is the initial tangent vector of the unique geodesic ray which starts at $\widetilde \alpha(t)$ and has limit point $\gamma(\infty) \in S_\infty$.   Then $v\alpha : I \to SM$ is defined to be the composite of $\widetilde{v\alpha}$ with the map on unit tangent bundles $S\HH \to SM$ given by the  universal cover $\HH \to M$.  It is not difficult to show that $v\alpha$ is independent of the choice of lift of $\widetilde \alpha$ and of the identification of the universal cover of $M$ with hyperbolic space.

Farrell and Jones' insight was to see how the asymptotic transfer behaves with respect to the geodesic flow $g^t : SM \times \R \to \R$.  Note that the geodesic flow is not an isometry, in fact it is an Anosov flow.  To quote \cite[p. 21]{FJ(1990)} ``The key property of the asymptotic transfer is that the geodesic flow shrinks $v\alpha$ in every direction except the flow line direction; i.e.~it deforms arbitrarily close to a flow line (as $t \to +\infty$) while keeping it bounded above in length.''

Here is a vague outline of how Farrell and Jones studied the Whitehead group of a fundamental group of a hyperbolic manifold and eventually showed it is zero.  Take an $h$-cobordism with base $M$.  Pull it back to an $h$-cobordism over $SM$.  Give $SM$ the metric induced by pulling back the standard metric by $g^t$ for $t$ large.   Given a deformation retract of an $h$-cobordism of $M$ with tracks $\alpha$, it is not difficult to construct a deformation retract of the pulled back $h$-cobordism with tracks the asymptotic transfer $v\alpha$.   Then the tracks are not small, but rather close to the flow lines.  These flow lines are either a one-to-one immersed image of the real numbers or diffeomorphic to the circle.  In either case the Whitehead group has been already shown to be trivial.  In this way, Farrell and Jones show that the pulled back $h$-cobordism is trivial.  In the case the manifold is odd-dimensional, the unit spheres in the tangent space have Euler characteristic 2, and Farrell and Jones concluded that $2 \wh(\pi_1M) = 0$.  However, to show that the Whitehead group vanishes Farrell and Jones  replaced the sphere bundle by a  bundle over the (noncompact) space $M\times \R$ whose fiber has Euler characteristic one.

The technical difficulties that Farrell and Jones had to overcome to make this outline a proof cannot be underestimated.  They had to formulate and prove sophisticated foliated control theorems \cite{FJ(1986)} adapted to the flow line situation.  They had to construct this noncompact bundle, the {\em enlargement} of $M$.  And, in the $L$-theory situation, one needs a fiber with signature one, and Farrell and Jones were forced to use $SM \times_{\Z_2} SM$ which is an orbifold, not a manifold.

Farrell and Jones continued their study of topological rigidity.  Using the methods outlined above, they showed the Borel Conjecture holds for negatively curved manifolds and applied their results to locally symmetric spaces.  They extended their results to the nonpositively curved case \cite{FJ(1993)rigid}, using a new idea, that of focal point transfer.  Instead of transferring using a point at infinity, they transfer using a point in $\widetilde M$ which is far away from the curve $\alpha$.  And, in their last joint work on topological rigidity \cite{FJ(1998)}, they used Cheeger-Fukaya-Gromov collapsing theory to study topological rigidity for noncompact manifolds, deducing, among other things, that the Borel Conjecture holds for complete closed affine manifolds of dimension greater than four.  Another nice application in this paper is to a Nielsen realization question.  Farrell and Jones show that if $M^n$ is a closed negatively curved manifold of dimension greater than four and if $F$ is a finite subgroup of Out$(\pi_1M)$ whose inverse image in Aut$(\pi_1M)$ is torsion-free, then $F$ acts on $M$ inducing the given outer automorphism of the fundamental group.  The Nielsen realization question without the torsion-free hypothesis seems quite interesting, but difficult given the example constructed in \cite{FL(2004)} by Farrell and Jean Lafont and an example of Jonathan Block and Shmuel Weinberger which shows Nielsen realization does not hold in general.

The next topic I will mention here is the stable pseudoisotopy space.  The {\em pseudoisotopy space} of a compact manifold $M$ is $$P(M) = \text{Homeo}(M \times [0,1] \text{ rel } M \times 0)$$ and the {\em stable pseudoisotopy space} is $${\cal P}(M) = \colim_{n \to \infty} P(M \times I^n).$$  The stable pseudoisotopy space is related to homeomorphism and diffeomorphism groups as well as to Waldhausen $A$-theory.  The case $M = S^1$ seem particularly interesting.  Waldhausen showed the homotopy groups of ${\cal P}(M)$ are rationally trivial and Kiyoshi Igusa showed that $\pi_0{\cal P}(S^1) \cong \oplus_\infty \Z/2$.  Farrell and Jones \cite{FJ(1987)} proved that for a closed negatively curved manifold, ${\cal P}(M) = \colim_{n\to \infty} {\cal P}(S^1)^n$ by using the asymptotic transfer again to gain control at the flow lines of the geodesic flow, but with the difference that the closed geodesics in $M$ each contribute a factor of ${\cal P}(S^1)$.

\section{Topology in service of geometry}

Let $\Sigma^n$ be an exotic sphere and $T^n$ the flat $n$-torus.  Let  $n \geq 5$.  Surgery theory shows that the manifolds $T^n$ and $T^n \#  \Sigma^n$ are not diffeomorphic.  This uses three ingredients:  every self-homeomorphism of the torus is homotopic to a diffeomorphism, the Borel Conjecture for $T^n \times I$, and the (stable) parallelizability of the torus\footnote{The Borel Conjecture for a compact manifold $N$ with boundary asserts that any homotopy equivalence from a compact manifold to $N$ which is  already a homeomorphism on the boundary is homotopic, relative to the boundary, to a homeomorphism.  The first two ingredients above show that any self-homeomorphism of the torus is topologically psuedoisotopic to a diffeomorphism, then smoothing theory gives the nondiffeomorphism result.}.
Furthermore, for every $\varepsilon > 0$, one can place a Riemannian metric on $T^n \#  \Sigma^n$ and rescale so that the sectional curvatures lies in the interval $[-\varepsilon ,\varepsilon ]$ (as can be for any Riemannian manifold).    A more interesting geometric question would be to ask for curvature pinched near zero, while keeping the diameter bounded, that is, to put an ``almost flat'' structure on $T^n \# \Sigma^n$.  However Michael Gromov and Ernst Ruh proved that almost flat manifolds are infranilmanifolds.  Any infranilmanifold with abelian fundamental group is diffeomorphic to a torus.  Hence $T^n \#  \Sigma^n$ does not admit a metric with bounded diameter and pinched curvature.

Farrell and Jones \cite{FJ(1989)exotic} studied  analogous questions for closed hyperbolic manifolds $M^n$; however this was much more difficult, important, and surprising.  Sullivan, building on joint work with Pierre Deligne,  showed that any hyperbolic manifold has a stably-parallelizable finite cover $\widehat M \to M$.  Furthermore, any self-homeomorphism of $\widehat M$ is homotopic to a diffeomorphism (Mostow rigidity) and Farrell and Jones had already proved topological rigidity for $\widehat M \times I$.   It follows that   $\widehat M \#  \Sigma$ and $\widehat M$ are not diffeomorphic.  However, the metric poses much more difficult problems.  Through an intricate geometric construction, they showed that for any closed hyperbolic $M$, for any $\varepsilon  > 0$, there exists a finite cover $\widehat M \to M$ so that $\widehat M \# \Sigma$ is not diffeomorphic to $\widehat M$ and a metric on $\widehat M \# \Sigma$ whose sectional curvatures lie in the interval $[-1-\varepsilon ,-1+\varepsilon ]$.  This is startling for a number of reasons.  First, note that the manifolds are obviously homeomorphic.   Second,  by Mostow Rigidity, $\widehat M \# \Sigma$ cannot admit a metric of constant negative curvature, or else it would be isometric, hence diffeomorphic to $\widehat M$.  Thus these manifolds can be added to the short list of closed manifolds which have (pinched) negative curvature and are not diffeomorphic to a locally symmetric space.  Third, Blaine Lawson and S.~T.~Yau conjectured that two homotopy equivalent negatively curved closed  manifolds are diffeomorphic.  Thus the examples of Farrell and Jones provide counterexamples to the Lawson-Yau conjecture.

I digress a bit and explain why the Lawson-Yau conjecture was well-founded and thus why a counter-example to it was so surprising.  The conjecture was motivated by the study of harmonic maps.  Work of  Jim Eells, Joseph Sampson, S.~I.~Al'ber, and Philip Hartmann showed that any map $f : M \to N$ between negatively curved closed manifolds has a unique harmonic representative $f_h : M \to N$ in its homotopy class, provided $f_*(\pi_1M)$ is not cyclic.  Thus if $f : M \to N$ is a homotopy equivalence with homotopy inverse $g: N \to M$, it is natural to think that $f_h$ and $g_h$ are mutual inverses, especially since Farrell and Jones proved that $f$ and $g$ are homotopic to homeomorphisms \cite{FJ(1993)rigid}.   In the special case of hyperbolic manifolds, $f_h$ and $g_h$ are inverses by Mostow rigidity, since isometries are of course harmonic.  But, to quote Farrell who paraphrased Gottfried Leibniz, ``the best of all possible maps is sometimes not good enough."

Farrell and Jones (and collaborator Pedro Ontaneda) wrote a sequence of deep papers exploring this phenomena, examining, for example, the extent to which $f_h$ and $g_h$ are homeomorphisms (in  \cite{FJO(1998)} they discussed examples with $f_h$ not a homeomorphism and in \cite{FOR(2000)} examples where, in addition, $f$ was a diffeomorphism) and generalized from hyperbolic manifolds to other symmetric spaces.

The final topic I mention in this section is the study of the space of negatively curved metrics and geometries.  Let $M^n$ be a closed smooth manifold and let $\met(M^n)$ be the space of all Riemannian metrics on $M^n$.  Then since any two metrics $g_0$ and $g_1$ are connected by a line segment $tg_0 + (1-t)g_1$, the space $\met(M^n)$ is contractible.  Let $\met^{\sectional < 0 }(M^n)$ be the subspace consisting of all negatively curved metrics.  Richard Hamilton's theorem on Ricci Flow shows that for $n = 2$, $\met^{\sectional < 0 }(M^n)$ is contractible.  Farrell and Ontaneda \cite{FO(2010)} show that $\met^{\sectional < 0 }(M^n)$ is not path connected for all negatively curved manifolds of dimension greater than nine and in \cite{FO(2009)} they exhibit the existence of high-dimensional hyperbolic manifolds where the space ${\cal T}^{\sectional < 0 }(M^n)$ of all (marked) negatively curved geometries on $M^n$ is not contractible; ${\cal T}^{\sectional < 0 }(M^n)$ is the quotient space of $\met^{\sectional < 0 }(M^n)$ obtained by identifying isometric metrics (by isometries homotopic to the identity).   They are pursuing these ideas in a sequence of preprints.  A notable feature of this recent work is the application of pseudoisotopy results of Hatcher, Waldhausen, and Igusa.

\section{Farrell-Jones Isomorphism Conjecture}

The isomorphism conjecture of Farrell and Jones \cite{FJ(1993)isom} marks a return to the theme of geometry in service of topology.  The story is amazing -- the geometry of negatively curved manifolds and their geodesics led Farrell and Jones to a general conjecture for computing $K_*(\Z\pi)$ and $L_*(\Z\pi)$ for an arbitrary group $\pi$.  Recall these groups are the key ingredients for classifying manifolds with fundamental group $\pi$.  I need to set the stage again, reminding the reader what is known for torsion-free groups, finite groups, and amalgamated products and HNN extensions.  For fundamental groups of aspherical manifolds and more generally torsion-free groups, a strong version of the Borel Conjecture predicts homological behavior, namely that assembly maps
\begin{align*}
H_*(B\pi: {\bf K}(\Z)) & \xrightarrow{\cong}  K_*(\Z\pi)\\
H_*(B\pi: {\bf L}(\Z)) & \xrightarrow{\cong}  L_*(\Z\pi)
\end{align*}
are isomorphisms. However, for finite groups the computations of $K_*(\Z\pi)$ and $L_*(\Z\pi)$ depend on representation theory and number theory, not on homology.    The theorems of Waldhausen and Cappell pointed to another sort of non-homological behavior, namely nilpotent phenomena for amalgamated products and HNN extensions.

It was hard to imagine that there could be a general picture.  But Farrell and Jones let geometry lead the way.  The fundamental group of a closed hyperbolic manifold is determined by a discrete torsion-free cocompact subgroup of the isometry group $\isom(\HH^n)$.  Farrell and Jones studied how the asymptotic transfer behaves for an arbitrary discrete cocompact subgroup $G$ of $\isom(\HH^n)$.  As in the torsion-free case, the geodesic flow in the unit sphere bundle $S(\HH^n)$ concentrates arcs in the image of the asymptotic transfer near the flow lines of the geodesic flow.  But the group theory is more complex.  One can show that a subgroup of $G$ which leaves a geodesic invariant must be virtually cyclic, and one can show that any virtually cyclic group must be finite, a semidirect product $F \rtimes \Z$ with $F$ finite, or an amalgamated product $G_0 *_F G_1$ with $F$ finite and of index two in $G_0$ and $G_1$.  Thus these are precisely the (sub)groups for which their controlled topology methods do not give any information.  The isomorphism conjecture says that $K$ and $L$-theory behave homologically except for non-homological behavior contributed by virtually cyclic subgroups.  In the modern formulation (due to myself and Wolfgang L\"uck), the Farrell-Jones Conjecture says that for an arbitrary group $G$, the assembly maps
\begin{align*}
H_*^G(E_{\vc}G; {\bf K}) \xrightarrow{\cong} K_*(\Z G)\\
H_*^G(E_{\vc}G; {\bf L}) \xrightarrow{\cong} L_*(\Z G)
\end{align*}
are isomorphisms.  Here $E_{\vc}G$ is the classifying space for $G$-actions with virtually cyclic isotropy.
This gives a conjectural view on what the classification of manifolds looks like, although computing the left hand sides is non-trivial because of Nil and UNil issues and because one must compute this version of group homology.   And, even once one computes the algebraic $K$- and $L$-groups, there are other surgery theoretic issues which arise when classifying a particular manifold.  L\"uck and I carried out this extended program for certain torus bundle over lens spaces.

Farrell and Jones \cite{FJ(1993)isom} proved a pseudoisotopy version of the isomorphism conjecture for discrete cocompact virtually torsion-free subgroups of the isometry group of the universal cover of a closed non-positively negatively curved manifold.  Recently L\"uck and collaborators Arthur Bartels and Holger Reich have proved the Farrell-Jones conjecture in K-theory for hyperbolic groups, in L-theory for both hyperbolic and
CAT(0) groups, and a version of it in lower K-theory for CAT(0) groups.

\section{Epilogue}

I first met Tom Farrell and Lowell Jones in 1989 at a diner in Poland where they shared a small table with my wife and my two daughters under the age of five.  We did not discuss topological rigidity at that lunch.
Although, I have only met Lowell a handful of times,
Tom has become a great friend and a mathematical inspiration.  I am honored to have helped co-organize the conference in Morelia in their honor and to co-edit this volume.

\bigskip
\noindent
Jim Davis

\bigskip

\noindent James F. Davis\\
Department of Mathematics\\ Indiana University \\Bloomington,
Indiana 47405 \\ USA \\Email: jfdavis@indiana.edu

\end{document}